\documentclass[a4paper,fleqn]{cas-sc}

\usepackage[numbers]{natbib}

\usepackage{graphicx}
\usepackage{subcaption}
\usepackage{parcolumns}
\def\tsc#1{\csdef{#1}{\textsc{\lowercase{#1}}\xspace}}
\tsc{WGM}
\tsc{QE}

\begin{document}
\let\WriteBookmarks\relax
\def\floatpagepagefraction{1}
\def\textpagefraction{.001}

\shorttitle{High-Order Interior Penalty Finite Element Methods for Fourth-Order
Phase-Field Models in Fracture Analysis}    

\shortauthors{Tian Tian, Chen Chunyu, Wei Huayi}  

\title [mode = title]{High-Order Interior Penalty Finite Element Methods for Fourth-Order
Phase-Field Models in Fracture Analysis}  



%

\author[1]{Tian Tian}[orcid=0009-0005-3440-9220]
\fnmark[1]
\ead{tiantian@smail.xtu.edu.cn}
\affiliation[1]{organization={School of Mathematics and Statistics, Xiangtan University},
                postcode={411105}, 
                city={Xiangtan},
                state={Hunan},
                country={China}}
\fntext[1]{The author was supported by the Graduate Innovation
Project of Xiangtan University (No. XDCX2023Y135).}

\author[1]{Chen Chunyu}
\ead{cbtxs@smail.xtu.edu.cn}

\author[2]{Huayi Wei}
\fnmark[2]
\ead{weihuyai@xtu.edu.cn}
\cormark[1]
\affiliation[2]{organization={School of Mathematics and Statistics, Xiangtan
    University; National Center of Applied Mathematics in Hunan; Hunan Key Laboratory for
Computation and Simulation in Science and Engineering},
                city={Xiangtan},
                postcode={411105},
                state={Hunan},
                country={China}}
\cortext[1]{Corresponding author}

\fntext[2]{Science Foundation of China (NSFC) (Grant Nos. 12371410, 12261131501) and the
construction of innovative provinces in Hunan Province (Grant No. 2021GK1010).}


\begin{abstract}
    This paper presents a novel approach for solving fourth-order phase-field
    models in brittle fracture mechanics using the Interior Penalty Finite
    Element Method (IP-FEM). The fourth-order model improves numerical stability
    and accuracy compared to traditional second-order phase-field models,
    particularly when simulating complex crack paths. The IP-FEM provides an
    efficient framework for discretizing these models, effectively handling
    nonconforming trial functions and complex boundary conditions.  

    In this study, we leverage the FEALPy framework to implement a flexible
    computational tool that supports high-order IP-FEM discretizations. Our
    results show that as the polynomial order increases, the mesh dependence of
    the phase-field model decreases, offering improved accuracy and faster
    convergence. Additionally, we explore the trade-offs between computational
    cost and accuracy with varying polynomial orders and mesh sizes. The
    findings offer valuable insights for optimizing numerical simulations of
    brittle fracture in practical engineering applications.
\end{abstract}



\maketitle

\section{Introduction}
In recent years, the phase-field method has emerged as a powerful tool for
simulating fracture processes\cite{Wu2018}. Among the different approaches, fourth-order
phase-field models have gained attention for their ability to improve numerical
stability and accuracy, especially when simulating complex crack paths.
Compared to traditional second-order phase-field models\cite{Wu2020, Nguyen2015,
Wu2020book}, the fourth-order phase-field model incorporates higher-order
derivative terms to enhance numerical stability and accuracy. This improvement
is particularly evident when simulating complex crack paths, where the model
minimizes nonphysical oscillations, thereby increasing computational reliability
and efficiency \cite{Borden2014, Wu2017}. By contrast, second-order phase-field
models often struggle with numerical instability, strong mesh dependency, and
high computational costs during the simulation of intricate crack propagation
\cite{Tanne2018}.

To address these challenges, significant progress has been made in advancing
fourth-order phase-field modeling. For instance, Amiri et al. \cite{Amiri2016}
proposed a fourth-order phase-field model leveraging the local maximum entropy
(LME) approximation. This approach directly solves fourth-order governing
equations by constructing high-order continuity functions, eliminating the need
for traditional double second-order decomposition and enabling efficient crack
path resolution on coarse meshs. Similarly, Borden et al. \cite{Borden2014}
developed thermodynamically consistent governing equations based on variational
principles. This approach improved solution smoothness, accelerated numerical
convergence, and showcased the potential of fourth-order models in
three-dimensional fracture problems. Moreover, recent researchers have adopted
hybrid solving strategies by integrating continuous and discontinuous Galerkin
methods \cite{Svolos2022}, significantly enhancing numerical stability for
dynamic fracture problems. Additionally, adaptive meshless algorithms combined
with mesh refinement techniques have emerged as promising tools for reducing
computational costs while maintaining high accuracy \cite{Shao2024}.

Efficient numerical solutions for fourth-order phase-field models require
appropriate computational techniques\cite{Hughes2012, Babuska1973}. In this regard, the Interior Penalty
Finite Element Method (IP-FEM) provides a flexible and efficient
framework\cite{Scott1990}.
Originally proposed by Douglas and Dupont \cite{Douglas1976}, the IP-FEM
addresses the challenges of nonconforming trial functions in high-order elliptic
equations by introducing penalty terms. This approach has proven to be a robust
solution for problems involving complex boundaries and interfaces. Subsequent
advancements by Wheeler \cite{Wheeler1978} provided in-depth analyses of
the effects of penalty parameters on solution stability and convergence. Arnold
\cite{Arnold1982} further unified the discontinuous Galerkin framework
and demonstrated the superior performance of IP-FEM in handling high-order
derivative problems.

This work explores the application of IP-FEM in solving fourth-order phase-field
models by leveraging the capabilities of the FEALPy software package
\cite{Wei2017}. A general programmatic framework is implemented to support
arbitrary high-order IP-FEM discretizations, offering a robust and efficient
tool for solving the fourth-order phase-field fracture model. Through rigorous
numerical experiments, we validate IP-FEM's superior efficiency and accuracy in
addressing complex fracture problems. Notably, we observe that as the finite
element degree increases, the mesh dependence of fourth-order phase-field models
is further reduced, highlighting their advantages. Finally, a comprehensive
comparison and analysis of numerical performance under varying polynomial orders
and mesh resolutions are conducted. This study provides theoretical insights and
practical guidelines for achieving low-cost, high-accuracy simulations of
complex crack propagation paths, offering a robust and efficient approach for
simulating complex crack propagation in engineering applications.

The remainder of this paper is organized as follows: In Section \ref{sec:model}, we present
the theoretical formulation of the fourth-order phase-field model and the
governing equations for crack propagation. Section \ref{sec:algorithm} is dedicated to the
numerical discretization using the nonlinear Interior Penalty Finite Element
Method (IP-FEM). Section \ref{sec:experiments} discusses the numerical
experiments, including the problem setups, boundary conditions, and results from
various simulations.  Finally, Section \ref{sec:conclusion} concludes the paper and outlines
potential directions for future research.

\section{Mathematical model}
\label{sec:model}
\subsection{Hybrid model for crack propagation}
The phase-field model employs a continuous field variable \( d(x) \) to
represent cracks within a material. 
Borden et al. \cite{Borden2012} proposed a
fourth-order phase-field model for the crack surface density, expressed as:

\[
\gamma(d, \nabla d, \Delta d) = \frac{1}{4l_0} d^2 + \frac{l_0}{2} |\nabla d|^2 + \frac{l_0^3}{32} (\nabla^2 d : \nabla^2 d)
\]

Here, \( l_0 \) is a scale factor that controls the width of the crack. 
We employ the Hybrid model \cite{Ambati2015} for the positive and
negative decomposition of strain energy. In this model, we
define\cite{Bourdin2000}:
\begin{equation*}
e_s^+(\boldsymbol{\varepsilon}) = e_s(\boldsymbol{\varepsilon}) = \frac{\lambda}{2} \, \text{tr}(\boldsymbol{\varepsilon})^2 + \mu \, \boldsymbol{\varepsilon} : \boldsymbol{\varepsilon},
\quad e_s^-(\boldsymbol{\varepsilon}) = 0.
\end{equation*}
where \( \lambda \) is Lamé's first parameter, \( \mu \) is Lamé's second
parameter (the shear modulus).

To prevent reversible cracking, for the positive strain energy in the
phase-field equation, we use the maximum history strain field function proposed
by Miehe \cite{Miehe2016}:
\begin{equation}
\mathcal{H}(\boldsymbol{x}, t) = \max_{s \in [0,t]} e_s^+(\boldsymbol{\varepsilon}(\boldsymbol{x}, s)).
\end{equation}

Here, \( e_s^+(\boldsymbol{\varepsilon}) \) is expressed as:
$e_s^+(\boldsymbol{\varepsilon}) = \frac{\lambda}{2} \langle
\text{tr}(\boldsymbol{\varepsilon}) \rangle_+^2 + \mu \,
\text{tr}(\boldsymbol{\varepsilon}_+^2)$, and
$\boldsymbol{\varepsilon}_{\pm} = \sum_{i=0}^{n-1} \langle \varepsilon_i
\rangle_{\pm} \, \boldsymbol{n}_i \otimes \boldsymbol{n}_i$. Where, \(
\varepsilon_i \) represents eigenvalues, \( \boldsymbol{n}_i \otimes
\boldsymbol{n}_i \) represents the eigenvectors, \( \langle \cdot \rangle_{\pm}
\) denotes the Macaulay bracket, defined as: $\langle x \rangle_{\pm} =
\frac{1}{2} \left( x \pm |x| \right).$

The resulting governing equations are:
\begin{equation}
\begin{aligned}
\begin{cases}
\rho \ddot{\boldsymbol{u}} - \text{div} \, \boldsymbol{\sigma} = \boldsymbol f, \\
2(1 - d) \mathcal{H} - G_c \left[ \frac{1}{2l_0} d - l_0 \Delta d + \frac{l_0^3}{16} \Delta^2 d \right] = 0.
\end{cases}
\end{aligned}
\end{equation}
with boundary conditions:
\begin{equation*}
\begin{cases}
\boldsymbol \sigma \cdot \boldsymbol n = \boldsymbol g , \quad \text{on } \partial \Omega_0, \\
d \cdot \boldsymbol n = 0, \quad \text{on } \partial \Omega.
\end{cases}
\end{equation*}

Here, $G_c$ is critical energy release rate, \( \boldsymbol{\sigma}(\boldsymbol{\varepsilon}) \) is the stress tensor, given by:
\begin{equation*}
\boldsymbol{\sigma}(\boldsymbol{\varepsilon}) = g(d) \frac{\partial e_s(\boldsymbol{\varepsilon})}{\partial \boldsymbol{\varepsilon}} = \left[ (1 - d)^2 + \epsilon \right] \left[ \lambda \, (\text{tr}(\boldsymbol{\varepsilon})) \, \boldsymbol{I} + 2\mu \, \boldsymbol{\varepsilon} \right],
\end{equation*}

\section{Algorithm design}
\label{sec:algorithm}

In the quasi-static crack model, the acceleration term $\ddot{\boldsymbol{u}}$ is neglected, yielding the variational form:
\begin{equation}
\begin{cases}
 (\boldsymbol{\sigma}(\boldsymbol{u}, d),
 \boldsymbol{\varepsilon}(\boldsymbol{v}))_{\Omega} = (\boldsymbol{f},
 \boldsymbol{v})_{\Omega} + \langle \boldsymbol{g}, \boldsymbol{v}
 \rangle_{\partial \Omega_N}, \\
 -2((1 - d) \mathcal{H}, \omega)_{\Omega} + \frac{G_c}{2l_0}(d, \omega)_{\Omega}
 + G_c l_0 (\nabla d, \nabla \omega)_{\Omega} 
 + \frac{G_c l_0^3}{16} (D^2 d, D^2
 \omega)_{\Omega} = 0.
\end{cases}
\end{equation}

Let $\Gamma_h$ denote the mesh set for the region $\Omega$ with mesh size $h$, and let $V_h$ be the Lagrange finite element space on $\Omega$ composed of polynomials of degree $p \geq 2$, i.e.,
\[
V_h = \{ v \in H^1(\Omega): v \in P_p(\tau), \forall \tau \in \Gamma_h \},
\qquad
V_h^0 = \{ v \in H^1_0(\Omega): v \in P_p(\tau), \forall \tau \in \Gamma_h \},
\]
where $P_p(\tau)$ represents the polynomial space of degree $p$ on each element $\tau$.

If $e \subset \Omega$, the unit normal vector
$\boldsymbol{n}_e$ is one of the two unit vectors normal to $e$, with the
direction from $\tau_-$ to $\tau_+$. On such an edge $e$, the following
definitions hold:
\[
\left[ \left[ \frac{\partial \boldsymbol{u}}{\partial \boldsymbol{n}} \right]
\right] = \frac{\partial \boldsymbol{u} \tau_+}{\partial \boldsymbol{n}}
\bigg|_e - \frac{\partial \boldsymbol{u} \tau_-}{\partial \boldsymbol{n}}
\bigg|_e,\qquad
\left\{ \frac{\partial^2 \boldsymbol{u}}{\partial \boldsymbol{n}^2} \right\} =
\frac{1}{2} \left( \frac{\partial^2 \boldsymbol{u} \tau_+}{\partial
\boldsymbol{n}^2} \bigg|_e + \frac{\partial^2 \boldsymbol{u} \tau_-}{\partial
\boldsymbol{n}^2} \bigg|_e \right).
\]

If $e \subset \partial \Omega$, the following holds:
\[
\left[ \left[ \frac{\partial \boldsymbol{u}}{\partial \boldsymbol{n}} \right] \right] = -\frac{\partial \boldsymbol{u}}{\partial \boldsymbol{n}_e}, \quad \left\{ \frac{\partial^2 \boldsymbol{u}}{\partial \boldsymbol{n}^2} \right\} = \frac{\partial^2 \boldsymbol{u}}{\partial \boldsymbol{n}_e^2}.
\]

The bilinear form $\mathcal{B}(\cdot, \cdot)$ is then defined as:
\begin{equation}
\begin{aligned}
\mathcal{B}(\boldsymbol{u}, \boldsymbol{v}) =  \sum_{\tau \in \Gamma_h}
\int_{\tau} D^2 \boldsymbol{u} : D^2 \boldsymbol{v} \, \mathrm{d} \boldsymbol{x} 
 + \sum_{e \in \mathcal{E}} \int_{e} ( \left\{ \frac{\partial^2
\boldsymbol{u}}{\partial \boldsymbol{n}^2} \right\} \left[ \left[ \frac{\partial
\boldsymbol{v}}{\partial \boldsymbol{n}} \right] \right] + \left[ \left[
\frac{\partial \boldsymbol{u}}{\partial \boldsymbol{n}} \right] \right] \left\{
\frac{\partial^2 \boldsymbol{v}}{\partial \boldsymbol{n}^2} \right\} +
\frac{\gamma}{|e|} \left[ \left[ \frac{\partial \boldsymbol{u}}{\partial
\boldsymbol{n}} \right] \right] \left[ \left[ \frac{\partial
\boldsymbol{v}}{\partial \boldsymbol{n}} \right] \right] ) \, \mathrm{d}
\boldsymbol{x}.
\end{aligned}
\end{equation}
Here, $\boldsymbol{n}$ is the unit normal, $\tau$ is an element, $e$ is an edge,
$\mathcal{E}$ is the set of edges, $\Gamma_h$ is the set of elements, and
$\gamma$ is the penalty parameter
\cite{Brenner2005}.

Then we can define the residuals as:
\begin{equation}
\begin{aligned}
\boldsymbol{R}_0 &= -\sum_{\tau \in \Gamma_h}
(\boldsymbol{\sigma}(\boldsymbol{u}_h, d_h),
\boldsymbol{\varepsilon}(\boldsymbol{v}_h))_{\tau} + \sum_{\tau \in \Gamma_h}
(\boldsymbol{f}, \boldsymbol{v}_h)_{\tau} 
+ \sum_{e \in \mathcal{E}_0} \langle
\boldsymbol{g}, \boldsymbol{v}_h \rangle_{e}, \quad \boldsymbol{v}_h \in
[V_h^0]^{\text{dim}}, \\ 
\boldsymbol{R}_1 &= \sum_{\tau \in \Gamma_h} 2((1 -
d_h) \mathcal{H}, \omega_h)_{\tau} - \sum_{\tau \in \Gamma_h}
\frac{G_c}{2l_0}(d, \omega_h)_{\tau} 
- G_c l_0 \sum_{\tau \in \Gamma_h} (\nabla
d_h, \nabla \omega_h)_{\tau} - \frac{G_c l_0^3}{16} \mathcal{B}(d_h, \omega_h),
\quad \omega_h \in V_h^0.
\end{aligned}
\end{equation}

Let $\{ \phi_m \}_{m=0}^{N}$ be the basis functions for the space $V_h$, where $N$ is the number of degrees of freedom. Let $\{ \Phi_{imj} \}_{m=0}^{N}$ be the basis functions for the space $[V_h]^{\text{dim}}$, where $i,j = 1,2,\dots,\text{dim}$, and
$
\Phi_{imj} = \delta_{ij} \phi_m,
$
where $\delta_{ij}$ is the Kronecker delta function.

Let $\boldsymbol{v}_h = \boldsymbol{\Phi}, \omega_h = \phi$. The residuals can then be expressed as:
\begin{equation}
\begin{aligned}
    \boldsymbol R_0 =& -\sum_{\tau \in
    \Gamma_h}(\boldsymbol{\sigma}(\boldsymbol{u}_h, d_h),
    \boldsymbol{\varepsilon}(\boldsymbol{\Phi}))_{\tau} + \sum_{\tau \in
    \Gamma_h}(\boldsymbol{f}, \boldsymbol{\Phi})_{\tau} 
    + \sum_{e \in
    \mathcal{E}_0} \langle \boldsymbol{g}, \boldsymbol{\Phi} \rangle_{e}, \\
    \boldsymbol R_1 =& \sum_{\tau \in \Gamma_h} 2 \left( (1 - d_h) \mathcal{H},
    \phi \right)_{\tau} - \sum_{\tau \in \Gamma_h} \frac{G_c}{2l_0} (d,
    \phi)_{\tau} 
    - G_c l_0 \sum_{\tau \in \Gamma_h} (\nabla d_h, \nabla
    \phi)_{\tau} - \frac{G_c l_0^3}{16} \mathcal{B}(d_h, \phi).
\end{aligned}
\end{equation}

The iteration is solved using the Newton-Raphson method:
\begin{equation}
\begin{aligned}
\begin{bmatrix}
-\frac{\partial \boldsymbol{R}_0}{\partial \boldsymbol{u_h}} & -\frac{\partial \boldsymbol{R}_0}{\partial d_h} \\
-\frac{\partial \boldsymbol{R}_1}{\partial \boldsymbol{u_h}} & -\frac{\partial \boldsymbol{R}_1}{\partial d_h}
\end{bmatrix}
\begin{bmatrix}
\Delta \boldsymbol{u} \\
\Delta d
\end{bmatrix} =
\begin{bmatrix}
\boldsymbol{R}_0 \\
\boldsymbol{R}_1
\end{bmatrix}
\end{aligned}
\end{equation}
where $\Delta \boldsymbol{u} = \boldsymbol{u_h}^{k} - \boldsymbol{u_h}^{k-1}$ and $\Delta d = d_h^k - d_h^{k-1}$. The components of the stiffness matrix are given by:
\begin{equation}
\left\{
\begin{aligned} 
    -\frac{\partial {\boldsymbol{R}_0}}{\partial \boldsymbol{u_h}}
    &= \sum_{\tau \in \Gamma_h} \frac{\partial \boldsymbol{\sigma}}{\partial
    \boldsymbol{\varepsilon}} (\nabla \boldsymbol{\Phi}^T, \nabla
    \boldsymbol{\Phi}), \\ 
    -\frac{\partial {\boldsymbol{R}_0}}{\partial d_h} &= \sum_{\tau \in
    \Gamma_h} \frac{\partial \boldsymbol{\sigma}}{\partial d_h} (\nabla
    \boldsymbol{\Phi}^T, \phi), \\ 
    -\frac{\partial \boldsymbol{R}_1}{\partial \boldsymbol{u_h}} &= \sum_{\tau
    \in \Gamma_h} -2(1 - d_h) \frac{\partial \mathcal{H}}{\partial
    \boldsymbol{u_h}} (\phi^T, \nabla \boldsymbol{\Phi}), \\ 
    -\frac{\partial \boldsymbol{R}_1}{\partial d_h} &= \sum_{\tau \in \Gamma_h}
    2 (\boldsymbol{\phi}^T \mathcal{H}, \phi) + \sum_{\tau \in \Gamma_h} G_c l_0
    (\nabla \phi^T, \nabla \phi)
    +\sum_{\tau \in \Gamma_h} \frac{G_c}{l_0}
    (\phi^T, \phi) + \frac{G_c l_0^3}{16}
\mathcal{B}(\phi, \phi).  \end{aligned}
\right.
\end{equation}

In a robust Staggered strategy, $\frac{\partial
\boldsymbol{R}_1}{\partial\boldsymbol{u_h}}$ and $\frac{\partial
\boldsymbol{R}_0}{\partial d_h}$ can be ignored, allowing the displacement
variable $\boldsymbol{u_h}$ and the phase field variable $d_h$ to be updated
independently. This independence simplifies the problem and allows for separate
updates of these variables without direct coupling in every iteration or time
step.

\section{Numerical experiments and results}
\label{sec:experiments}
In the numerical experiments, the penalty parameters were set as follows: for
order $p = 2$, $\gamma = 5$; for $p = 3$, $\gamma = 10$; and for $p = 4$,
$\gamma = 20$ \cite{Hughes2012, Arnold1982, Babuska1973}.

\begin{figure}[htbp]
    \centering
    \includegraphics[width=0.5\textwidth]{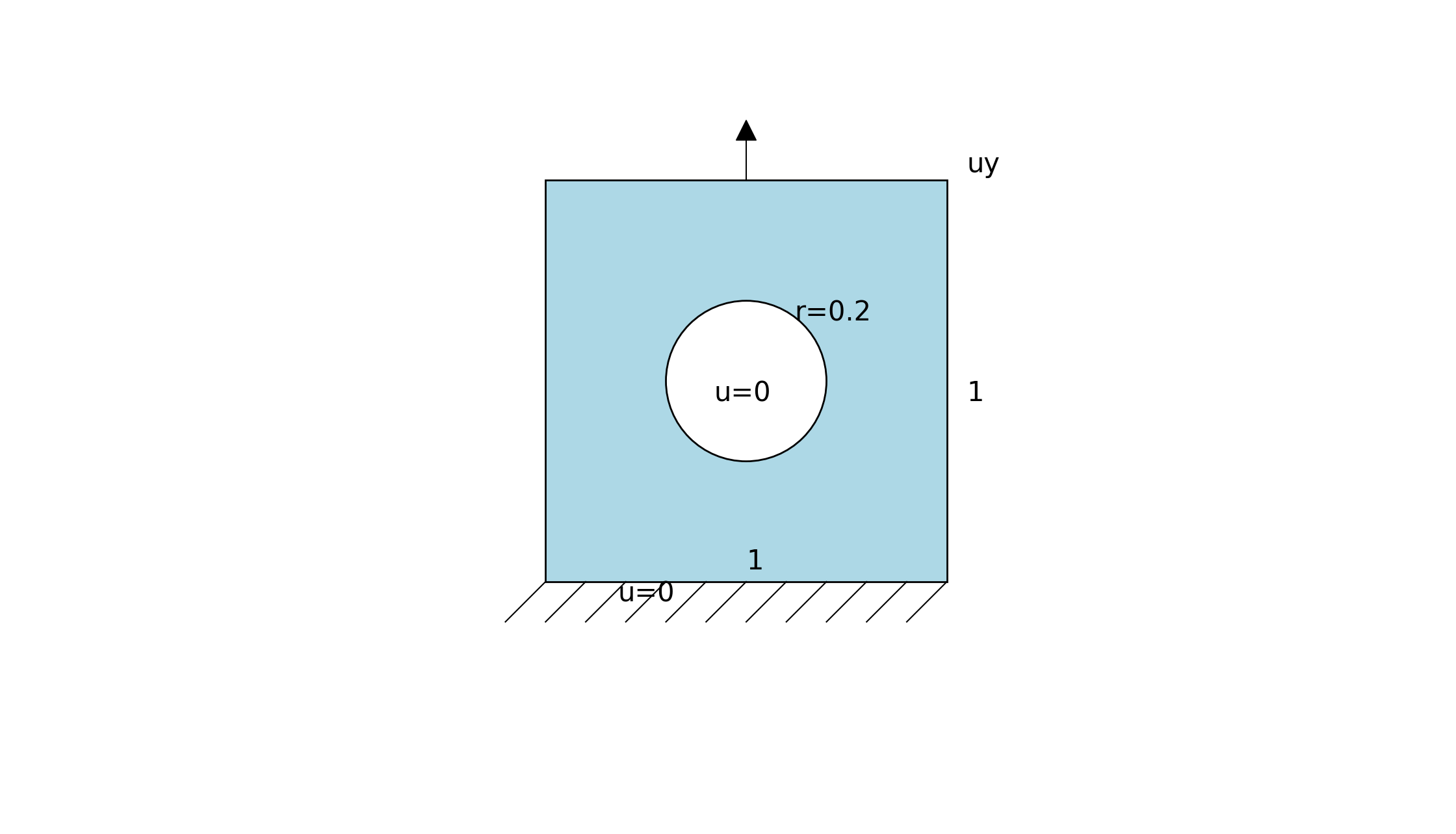}
    \caption{The square model with a circular hole (unit:cm)}
    \label{fig:ip_model0}
\end{figure}

This example examines a rigid circular inclusion within a square
plate, subjected to a vertical upward displacement applied to the top surface.
The domain \(\Omega\) is defined as the rectangular region \([0, 1] \times [0,
1]\), featuring a circular hole of radius \(0.2\) centered at the origin, as
depicted in Figure \ref{fig:ip_model0}. The material properties are specified as
follows: a critical energy release rate of \(G_c = 1\text{kN/mm}\), a length
scale factor of \(l_0 = 0.02\text{mm}\), Young's modulus of \(E =
200\text{kN/mm}^2\), and Poisson's ratio of \(\nu = 0.2\). The boundary
conditions include a Dirichlet boundary condition at the upper boundary (\(y =
1\)), where the displacement increases by \(\Delta u_y = 1.4 \times 10^{-2} \,
\text{mm}\) for the first 5 steps and then by \(\Delta u_y = 2.2 \times 10^{-3}
\, \text{mm}\) for the next 25 steps. Additionally, the displacement is set to
zero at the center of the circular notch.

\begin{figure}[htbp]
    \centering
    \begin{subfigure}{0.45\textwidth}
        \centering
        \includegraphics[width=\linewidth]{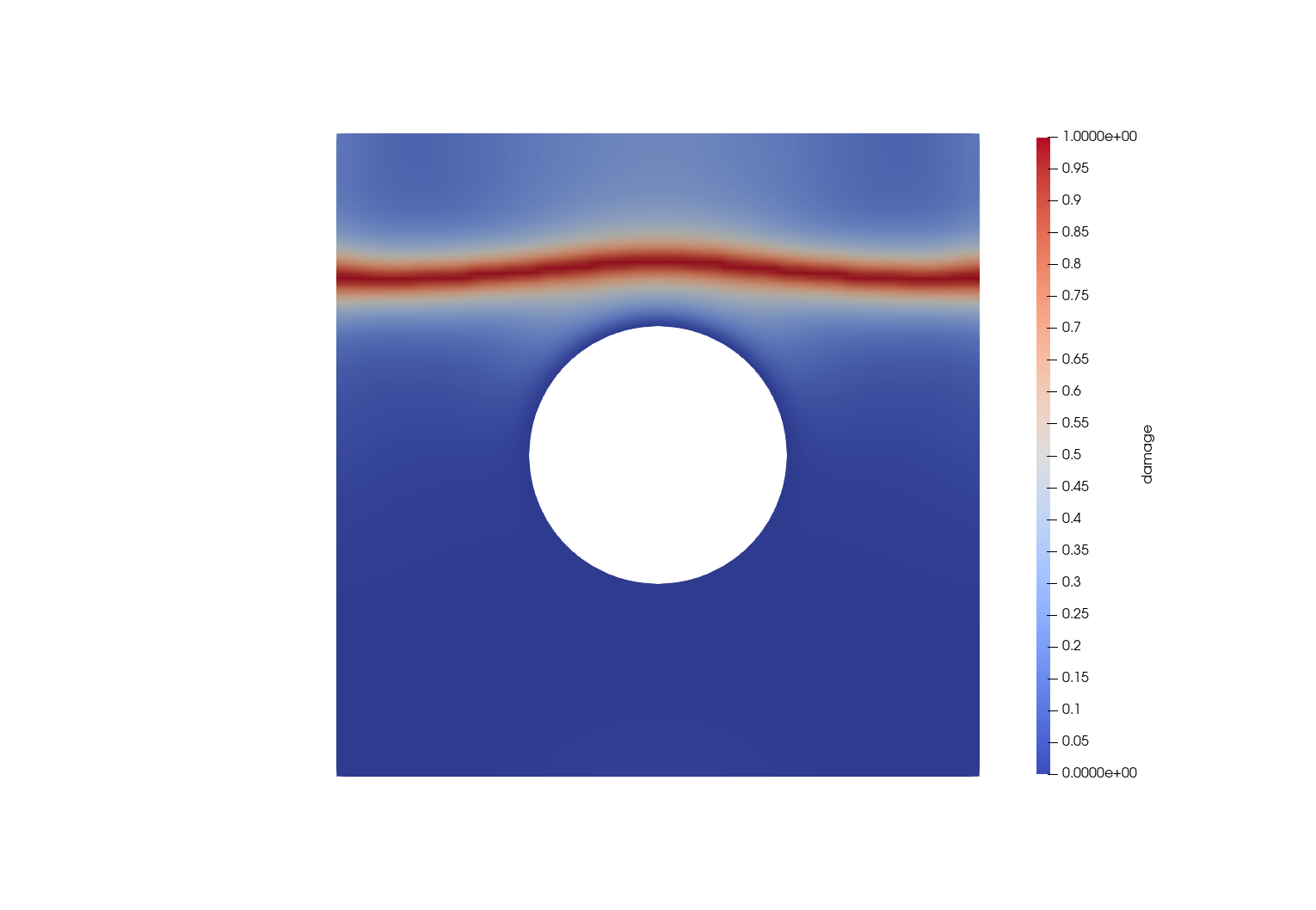}
    \end{subfigure}
    \begin{subfigure}{0.5\textwidth}
        \centering
        \includegraphics[width=\linewidth]{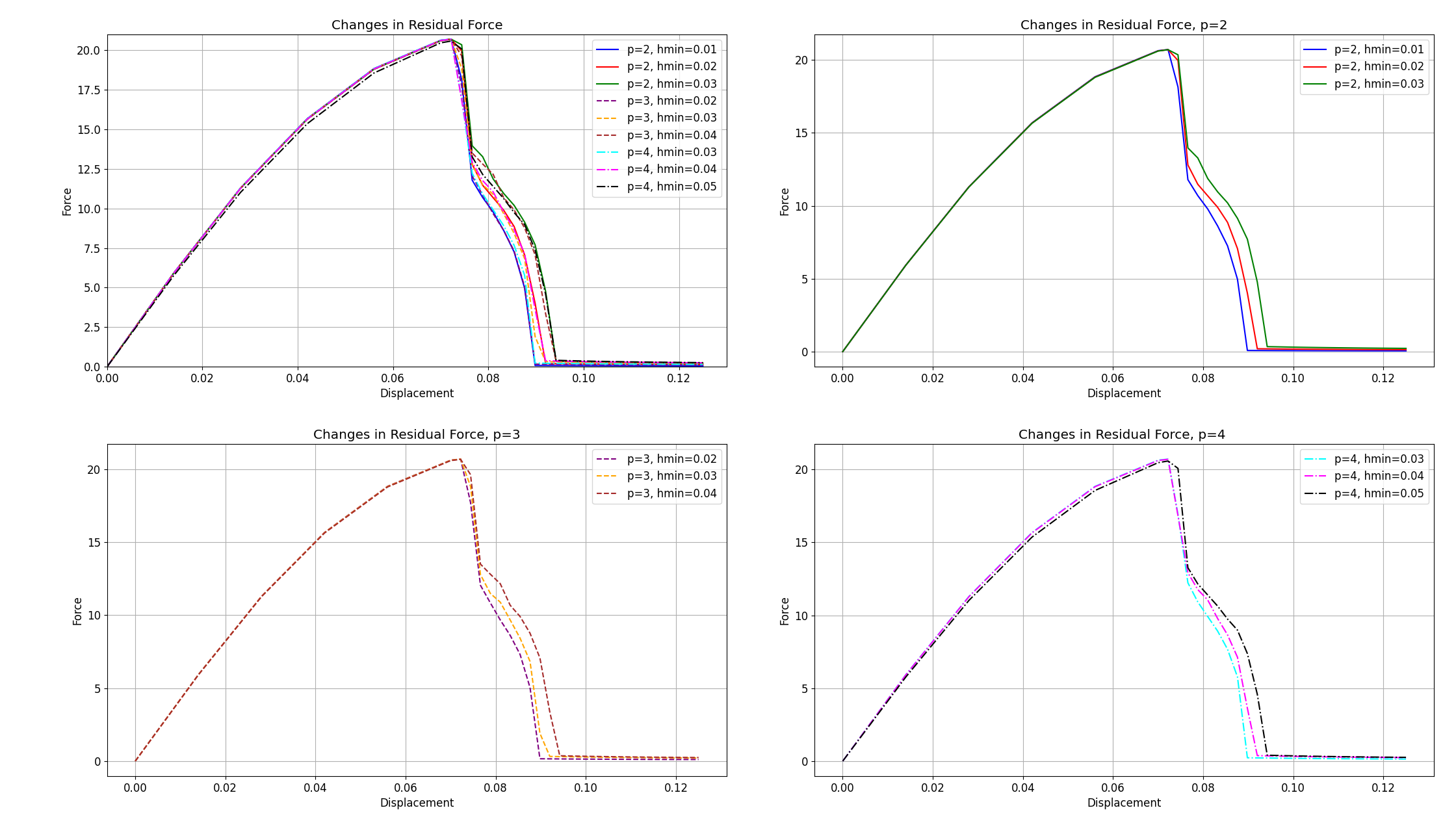}
    \end{subfigure}
    \caption{The final results (left) and the residual force for different
        degrees and mesh sizes in the square model with a circular hole
    (right).} 
    \label{fig:ip_model0_results}
\end{figure}

In this example, we investigate the performance and computational efficiency of
different finite element methods with varying polynomial degrees (\(p = 2, 3,
4\)) and mesh sizes (\(h_{\text{min}}\)) for simulating a phase-field fracture
model. Figure \ref{fig:ip_model0_results} shows the final fracture
morphology of the model.

As shown in Figure \ref{fig:ip_model0_results}, the residual force curves for
different polynomial degrees (\(p = 2, 3, 4\)) and mesh sizes (\(h = 0.01\) to
\(h = 0.05\)) reveal significant differences in terms of numerical stability and
convergence. 

For the lower-order finite element method (\(p = 2\)), the residual force curves
exhibit considerable variability, especially for finer meshes (\(h = 0.01\) to
\(h = 0.03\)). During the phase where the residual force rapidly decreases after
fracture, these curves show poor convergence, indicating that the lower-order
method struggles to provide accurate results without significantly finer mesh
sizes. This behavior is particularly evident in the early stages of crack
propagation, where the residual force curves for \(p = 2\) diverge more
noticeably compared to higher-order methods. This suggests that the lower-order
method is more sensitive to mesh refinement and requires finer meshes to achieve
stable and accurate results.

In contrast, the higher-order methods (\(p = 3\) and \(p = 4\)) demonstrate
superior numerical stability and convergence. The residual force curves for
these methods nearly coincide across different mesh sizes, indicating that
higher-order methods are less sensitive to mesh refinement. Notably, the \(p =
4\) method shows exceptional consistency across mesh sizes, suggesting that even
larger mesh sizes, such as \(h = 0.05\), can still maintain high accuracy and
good convergence. This is particularly evident in the post-fracture phase, where
the residual force curves for \(p = 4\) remain tightly clustered, regardless of
the mesh size. This behavior highlights the robustness of higher-order methods
in handling the complexities of crack propagation, even with coarser meshes.

Higher-order methods (\(p = 3\) and \(p = 4\)) significantly improve the
numerical accuracy and convergence of phase-field fracture simulations,
particularly for larger mesh sizes where lower-order methods (\(p = 2\))
struggle to maintain stability.  
The residual force curves for higher-order methods exhibit remarkable
consistency across different mesh sizes, demonstrating their robustness in
handling complex crack propagation problems.  

\section{Summary}
\label{sec:conclusion}
In this study, we investigated the numerical simulation of phase-field fracture
models using penalty finite element methods of varying orders. Our results
demonstrate that higher-order finite element methods (with $p=3$ and $p=4$) can achieve high accuracy even with coarser mesh sizes, making them a
highly efficient option for modeling fracture behavior in engineering
applications. Specifically, the analysis revealed that, for a given mesh size,
higher-order methods exhibited superior consistency and convergence properties
compared to the lower-order method ($p=2$), which required finer meshes to
achieve similar accuracy.

However, the results were sensitive to mesh refinement. Coarser meshes
introduced some numerical errors, especially in higher-order finite element
methods. While these errors remained within acceptable bounds for most practical
purposes, it is crucial to carefully select mesh sizes to balance computational
efficiency and accuracy. By optimizing mesh sizes and element orders, engineers
can achieve an optimal trade-off between solution accuracy and computational
cost.

Our future work will focus on extending this framework to more
complex geometries and dynamic fracture problems, where time-dependent behavior
plays a significant role. Additionally, the incorporation of material
nonlinearity and adaptive refinement strategies will be explored to further
enhance both the accuracy and efficiency of the simulations. These advancements
could significantly improve the versatility of the proposed method, enabling its
application to a wider range of real-world fracture problems in both engineering
and materials science.

\bibliographystyle{cas-model2-names}
\bibliography{references}  

\end{document}